\documentclass[mnsc,nonblindrev]{informs3}

\OneAndAHalfSpacedXII 

\renewcommand{\theARTICLETOP}{}

\usepackage{algorithm}
\usepackage{algorithmic}

\usepackage{natbib}
 \bibpunct[, ]{(}{)}{,}{a}{}{,}%
 \def\bibfont{\small}%
\TheoremsNumberedThrough     
\ECRepeatTheorems

\EquationsNumberedThrough    

\begin{document}


\RUNAUTHOR{D'Andreagiovanni et al.} 

\RUNTITLE{GUB Covers and Power-Indexed formulations for Wireless
Network Design}

\TITLE{GUB Covers and Power-Indexed formulations for Wireless
Network Design\thanks{
This is the authors' final version of the paper published in Management Science 59(1), 142-156, 2013.
DOI: 10.1287/mnsc.1120.1571 \; .
The final publication is available at INFORMS via http://pubsonline.informs.org/doi/abs/10.1287/mnsc.1120.1571}}

\ARTICLEAUTHORS{%
\AUTHOR{Fabio D'Andreagiovanni}
\AFF{Konrad-Zuse-Zentrum für Informationstechnik Berlin, 14195 Berlin, Germany, \EMAIL{d.andreagiovanni@zib.de}} 
\AUTHOR{Carlo Mannino}
\AFF{Department of Applied Mathematics, SINTEF, 0314 Oslo, Norway; and Department of Computer, Control, and
Management Engineering, Sapienza Università di Roma, 00185 Rome, Italy, \EMAIL{mannino@dis.uniroma1.it}}
\AUTHOR{Antonio Sassano}
\AFF{Department of Computer, Control, and Management Engineering, Sapienza Università di Roma, 00185 Rome, Italy, \EMAIL{sassano@dis.uniroma1.it}}
} 

\ABSTRACT{%
We propose a pure 0-1 formulation for the wireless network design
problem, i.e. the problem of configuring a set of transmitters to
provide service coverage to a set of receivers. In contrast with
classical mixed integer formulations, where power emissions are
represented by continuous variables, we consider only a finite set
of powers values. This has two major advantages: it better fits
the usual practice and eliminates the sources of numerical
problems which heavily affect continuous models.  A crucial
ingredient of our approach is an effective basic formulation for
the single knapsack problem representing the coverage condition of
a receiver. This formulation is based on the GUB cover
inequalities introduced by \cite{Wo90} and its core is an
extension of the exact formulation of the GUB knapsack polytope
with two GUB constraints. This special case corresponds to the
very common practical situation where only one major interferer is
present. We assess the effectiveness of our formulation by
comprehensive computational results over realistic instances of two typical technologies, namely WiMAX and DVB-T.
}%


 \HISTORY{Received August 1, 2011; accepted March 2, 2012, by Dimitris Bertsimas, optimization. Published
online in Articles in Advance August 20, 2012.}

\maketitle

%


\section{Introduction}\label{sec:introduction}


Wireless communication systems constitute one of the most
pervasive phenomena of everyday life. Television and radio
programs are distributed through broadcasting networks (both
terrestrial and satellite), mobile communication is ensured by
cellular networks, internet service is provided through broadband
access networks. Moreover, a number of security services are
provided by ad-hoc wireless networks. All these networks have
grown very rapidly during the last decades, generating  dramatic
congestion of radio resources such as frequency channels. Wireless networks provide
different services and rely on different technologies and
standards. Still, they share a common feature: they all need to
reach users scattered over an area with a radio signal that
must be strong enough to prevail against other unwanted interfering signals.

The perceived quality of service thus  depends on several signals,
wanted and unwanted, generated from a large number of transmitting
devices. Due to the increasing size of the new generation
networks, co-existing in an extremely congested radio spectrum and
subject to local and international constraints, establishing
suitable power emissions for all the transmitters has become a
very difficult task, which calls for sophisticated optimization
techniques.

Since the early 1980s several optimization models have been
developed to design wireless networks. It is claimed that the use
of automatic and optimization-oriented planning techniques may
lead to cost reduction of up to 30\% \citep{De05}. Concretely,
recent experiences have clearly shown that the adoption of
optimization techniques results in sensible increases in the
quality of coverage plans and in a more effective and efficient
use of the limited resources that a network administrator has at
disposal - see the case of \emph{atesio} for UMTS networks in
Germany \citep{atesio} and the case of the \emph{CORG} for DVB-T
networks in Italy \citep{CORG}.

Two fundamental issues must be faced when designing a wireless
network: localizing the transmitters and dimensioning their power
emissions. In most models, power emissions are represented as
continuous decision variables. This choice typically yields
ill-conditioned constraint matrices and requires the introduction
of very large coefficients to model disjunctive constraints. The
corresponding relaxations are very weak and state-of-the-art
Mixed-Integer Linear Programming solvers are often affected by
numerical instability. The use of continuous decision variables
also contrasts with the telecommunications practice. In fact, the
actual design specifications of real life antennas are always
expressed as rational numbers with bounded precision and,
consequently, assume a finite number of values.

Motivated by the above remarks we propose a pure 0-1 formulation
for the problem that is obtained by considering only a finite set
of power values. This formulation has two basic advantages: first,
the ensuing model better fits the usual practice and, second, the
numerical problems produced by the continuous variables are
sensibly reduced. Indeed, the new approach allows us to find
better solutions to large practical instances with less
computational effort. In addition, the model fits the common network
planning practice of considering a small number of power values
and it directly models power restrictions that are often imposed
by the technology (e.g., \citealt{MaDrHu07}). The situation
where only two power values (on, off) are allowed is not rare {\citep{Ri10}. Finally, the new approach easily allows for
generalizations of the model, such as power consumption
minimization or antenna diagram optimization.

For our purposes, a wireless network can be described as a set of
transmitters $B$ distributing a telecommunication service to a set
of receivers $T$.  A receiver is said to be {\em covered} (or {\em
served}) by the  network if it receives the service within a
minimum level of quality. The set $B$ actually contains all {\em candidate} transmitters: in general, only a subset of $B$ will  be activated to cover the set $T$. Transmitters and receivers are characterized by a number of locations and radio-electrical
parameters (e.g., geographical coordinates, power emission,
transmission frequency).   The {\em Wireless Network Design
Problem} (WND) consists of establishing suitable values for such
parameters with the goal of maximizing the coverage (or a revenue
associated with the coverage).

Each transmitter $b\in B$ emits a radio signal with power $p_b\in
[0,P_{\max}]$. We remark that a transmitter $b$ such that $p_b = 0$ is actually not activated and thus not deployed in the network. The power $p(t)$ received by receiver $t$ from
transmitter $b$ is proportional to the emitted power $p_b$ by a
factor $\tilde{a}_{tb} \in [0,1]$, i.e. $p(t) =
\tilde{a}_{tb}\cdot p_b$. The factor $\tilde{a}_{tb}$ is called
{\em fading coefficient} and summarizes the reduction in power
that a signal experiences while propagating from $b$ to $t$. The
value of a fading coefficient depends on many factors (e.g.,
distance between the communicating
devices, presence of obstacles, antenna patterns) and is commonly
computed through a suitable propagation model. For a detailed
presentation of all technical aspects, we refer the reader to
\cite{Ra01}.

To simplify the discussion, we assume here that all the
transmitters of the network operate at the same frequency. This
assumption is dropped in Section \ref{sec:results} where we introduce
the real-life application which motivated our developments. Among
the signals received from transmitters in $B$, receiver $t$ can
select a {\em reference signal} (or {\em server}), which is the
one carrying the service. All the other signals are interfering.

A receiver $t$ is regarded as served by the network, specifically
by server $\beta \in B$,  if the ratio of the serving power to the
sum of the interfering powers ({\em signal-to-interference ratio}
or \emph{SIR}) is above a threshold $\delta'$ \citep{Ra01}, the
\emph{SIR threshold}, whose value depends on the technology and
the desired quality of service:
\begin{equation}
\label{eq:firstSIRinequality} \frac{\tilde{a}_{t \beta} \cdot
p_{\beta}}{\mu + \sum_{b \in B \setminus\{\beta\}} \tilde{a}_{tb}
\cdot p_b} \geq \delta'.
\end{equation}

\medskip

\noindent
Note the presence of the system noise $\mu > 0$ among the
interfering signals. Since each
transmitter in $B$ is associated with a unique received signal,
in what follows we will also refer to $B$ as the set of signals
received by $t$. By letting
$\delta = - \mu \cdot \delta' < 0$ and letting:
\begin{eqnarray*}
  a_{tb} &=& \left\{
                    \begin{array}{ll}
                        \tilde{a}_{tb} \hspace{1.0cm}\mbox{ if $b = \beta$}
                        \\
                        \delta' \cdot \tilde{a}_{tb} \hspace{0.5cm} \mbox{ otherwise}
                    \end{array}
                \right.
\end{eqnarray*}

\noindent for every $b \in B$, inequality
(\ref{eq:firstSIRinequality})  can be transformed into the
so-called \emph{SIR inequality} by simple algebra operations:
\begin{equation}\label{eq:SIR}
 \sum_{b \in B\setminus\{\beta\}} a_{tb} \cdot p_b
 - a_{t \beta} \cdot p_{\beta} \leq \delta.
\end{equation}
For every $t\in T$, we have one inequality of type (\ref{eq:SIR}) for
each potential server $\beta\in B$. Receiver $t$ is served if
at least one of these inequalities is satisfied or, equivalently, if
the following disjunctive constraint is satisfied:
\begin{equation}\label{eq:disjunctive-SIR}
\bigvee_{\beta \in B} \left( \sum_{b \in
B\setminus\{\beta\}} a_{tb} \cdot p_b
 - a_{t \beta} \cdot p_{\beta} \leq \delta \right).
\end{equation}

\medskip

\noindent
The above disjunction  can be represented by a family of linear
constraints in the $p$ variables by introducing, for each $t\in T$
and each $b\in B$, a binary variable $x_{tb}$ that is equal to $1$
 if $t$ is served by $b$ and to $0$ otherwise. For each $\beta
\in B$, the following constraint is then introduced:
\begin{equation}\label{eq:SIR-BIGM}
 \sum_{b \in B\setminus\{\beta\}} a_{tb} \cdot p_b
 - a_{t \beta} \cdot p_{\beta}  - M  \cdot (1 - x_{t\beta})
 \leq \delta
\end{equation}

\noindent where $M$ is a large positive constant (big-\emph{M}).
When $x_{t\beta} = 1$ then (\ref{eq:SIR-BIGM}) reduces to
(\ref{eq:SIR}); when instead $x_{t\beta} = 0$  and $M$ is
sufficiently  large (for example, we can set $M = - \delta +
\sum_{b \in B\setminus\{\beta\}} a_{tb} \cdot P_{max}$),
(\ref{eq:SIR-BIGM}) is satisfied for any feasible power vector and
becomes redundant. Constraints of type (\ref{eq:SIR-BIGM}) appear
in the Mixed-Integer Linear Programs (MILP) for the WND presented
in several papers in different application contexts, such as radio
and video broadcasting (e.g., \citealt{MaRoSm06, MaMaSa09}), GSM
(e.g., \citealt{MaSc05}), UMTS (e.g., \citealt{AmEtAl06a, EiGe08,
KaKeOl06, NaEl10}), WiMAX \citep{Zh09}. Such MILPs are informally
called {\em big-M formulations}. For a comprehensive description
of the main elements that constitutes such models we refer to the
recent book by \cite{KeOlRa10} and to \cite{AmEtAl06b}. For a more
detailed discussion about how modeling an UMTS network, we refer
the reader to \cite{EiEtAl02}, and additionally to
\cite{SiEtAl06}, where focus is on dimensioning pilot channel
powers rather than the overall power emissions, considered as
fixed.

The MILPs have been also tailored to cope with uncertainty
affecting parameters of the model: in (\citealt{RoOl07}) and (\citealt{OlRo08}),
two stochastic optimization approaches are presented to establish
a robust location plan of the transmitters to tackle fluctuations
in the traffic demand; in (\citealt{HePr04}) and (\citealt{BiDa09}), Stochastic
and Robust Optimization are respectively adopted to tackle the
uncertainty affecting the fading coefficients.

WND instances of practical interest typically correspond to very
large MILPs. In principle, such programs can be solved by standard
Branch-and-Cut and by means of effective commercial solvers such
as IBM ILOG \cite{CPLEX}. However, it is well-known that the
presence of a great number of constraints of type
(\ref{eq:SIR-BIGM}) results in ill-conditioned instances, due to
the large variability of the fading coefficients, and weak bounds,
due to the presence of the big-{\em M} coefficients. Furthermore,
the resulting coverage plans are often unreliable (e.g.,
\citealt{KaKeOl06,KeOlRa10,MaMaSa09}). In some cases, feasible WND instances can be even considered as unfeasible.
 In practice, only small-sized WND instances can actually be solved to optimality.

It is interesting to note that though these problems are known,
only a limited number of papers of the wide literature about the
WND has tried to overcome them. \cite{KaKeOl06} proposed to
execute a post-processing procedure that tries to repair coverage
errors by eventually dropping service of a number of receivers.
\cite{NaEl10} focused on networks based on \emph{Code Division
Multiple Access (CDMA)} and adopted Benders' decomposition to
obtain a new problem where the big-\emph{M} coefficients are
eliminated. However, the fading coefficients are still present,
thus maintaining a relevant source of numerical problem. Inspired
by practical observations about DVB networks, \cite{MaMaSa09}
considered a relaxation of the WND problem, obtained by including
a single interfering transmitter in each SIR constraint and solved
by a heuristic approach. Finally, \cite{EiGe08} proposed a new
approach for reducing interference in an UMTS network to increase
the overall capacity of the network, under the assumption of
perfect power control.

All the previously cited work are based on modeling the power
emission of a transmitter as a continuous variable. In this paper,
we follow instead a different path: we discretize the continuous
power variables and consider only a finite number of feasible
values. We stress that discretization is a classical tool in
combinatorial optimization (e.g., \citealt{DyWo90}) and in
telecommunication modeling (e.g., \citealt{CaNoTr08,
FrWeDaKa08,MaDrHu07}), but, to our best knowledge, no effort has
been made to go beyond the simple use of discretized SIR
inequalities and replace them by more combinatorial inequalities.
By using discretization, we are instead able to completely
eliminate the two main sources of numerical issues, namely the
fading and the big-\emph{M} coefficients. We accomplish this by
introducing a set of (strong) valid inequalities for the resulting
0-1 problem that radically improve the quality of obtained
solutions. Additionally, solutions do not contain errors.


In the next section, we introduce our new contribution to the WND,
the Power-Indexed formulation. In Section \ref{sec:convex-hull}, we prove that for a special case that is very relevant in practice (single server interfered by a single transmitter), we can characterize the convex hull of the knapsack polytope associated with discrete power levels. In Section \ref{sec:algorithm}, we describe our
solution approach to the WND. Finally, extensive computational results on realistic
instances of WiMAX and DVB-T networks are presented in Section \ref{sec:results},
showing that the new approach outperforms the one based
on the big-\emph{M} formulation.

\section{A Power-Indexed formulation for the WND}\label{sec:power-indexed}

As discussed in the previous section, a classical and much
exploited model for the WND belongs to the class of the so-called
big-{\em M} formulations and writes as:
\begin{eqnarray}
    \makebox[13mm][l]{$\max $}
    \makebox[104mm][l]{$\displaystyle \sum_{t\in T}
        \sum_{b\in B}
        r_t \cdot x_{tb}$}
    && \hspace{3.5cm} (BM)
    \nonumber
    \\
    \makebox[13mm][l]{s.t.}
    \makebox[104mm][l]{$\displaystyle{\sum_{b \in B\setminus\{\beta\}} a_{tb} \cdot p_b
 - a_{t \beta} \cdot p_{\beta}  - M  \cdot (1 - x_{t\beta})
 \leq \delta}
    \hspace{2.0cm}
        t\in T, \hspace{0.1cm}
        \beta \in B$}
    \label{eq:BM-SIR}
    \\
    \makebox[104mm][l]{$\displaystyle{\sum_{b\in B} x_{tb} \leq
    1}
    \hspace{7.2cm}
        t\in T$}
    \label{eq:BM-oneserver}
    \\
    \makebox[104mm][l]{$\displaystyle{0 \leq p_b \leq P_{max}}
    \hspace{6.65cm}
    b \in B
    $}
    \nonumber
    \\
    \makebox[104mm][l]{$x_{tb} \in \{0,1\}
    \hspace{7.1cm}
    t \in T, \hspace{0.1cm} b \in B
    $}
    \nonumber
\end{eqnarray}

\noindent where $r_t$ is the revenue (e.g. population, number of
customers, expected traffic demand) associated with receiver $t
\in T$ and the objective function is to maximize the total revenue.
Constraint (\ref{eq:BM-SIR}) is the  SIR inequality
(\ref{eq:SIR-BIGM}) introduced in Section \ref{sec:introduction}
and constraint (\ref{eq:BM-oneserver}) ensures that each receiver
is served at most once.

Technology-dependent versions of (BM)  can be obtained from the
basic formulation by including suitable constraints or even new
variables. For example, in the case of WiMAX networks, a knapsack
constraint involving the service variables $x_{tb}$ is added to
(BM) to model the bandwidth capacity of each transmitter $b \in
B$ \citep{Zh09}.
In the case of antenna diagram
design, the number of power variables associated with each
transmitter $b$ is multiplied by 36 to represent the power
emissions along the 36 directions which approximate the horizontal
radiation pattern, and new constraints are
included to represent physical relations between different
directions (\citealt{MaMaSa09}).

As observed in the introduction, the model (BM) has serious
drawbacks both in terms of dimension of the solvable instances and
of numerical instability. We tackle these issues by
restricting the variables $p_b$ to assume value in the finite set
${\cal P} = \{P_1, \dots, P_{|{\cal P}|}\}$ of feasible power
values, with $P_1 = 0$ ({\em switched-off value}), $P_{|{\cal P}|}
= P_{max}$ and $P_i > P_{i-1}$, for $i = 2, \dots, |{\cal P}|$. To
this end, we introduce  a binary variable $z_{bl}$, which is 1 iff
$b$ emits at power $P_l$. Since $b$ is either switched-off or
emitting at a positive  value in ${\cal P}$, we have:
$$
\sum_{l\in L} z_{bl} = 1 \qquad b \in B
$$

\noindent
where $L = \{1,\dots,|{\cal P}|\}$ is the set of power
value indices or simply {\em power levels}. Then we can write:
\begin{equation}\label{eq:powervalue}
p_{b} = \sum_{l\in L} P_l \cdot z_{b l} \qquad b \in B.
\end{equation}
By substituting  (\ref{eq:powervalue}) in (\ref{eq:BM-SIR}), we
obtain  the following SIR constraint that only involves  0-1
variables:
$$
  \sum_{b  \in B\setminus\{\beta\}} a_{tb} \hspace{0.1cm} \sum_{l\in L} P_l \cdot z_{b l} \hspace{0.1cm}
- a_{t \beta} \hspace{0.1cm} \sum_{l\in L}  P_l \cdot z_{\beta l} \hspace{0.1cm}
- M
\cdot (1 - x_{t\beta})  \leq \delta
$$
The following \emph{discrete big-M formulation} $(DM)$ for
the WND with a finite number of power values directly derives from
(BM):

\begin{eqnarray}
    \makebox[13mm][l]{$\max $}
    \makebox[100mm][l]{$\displaystyle{\sum_{t\in T} \sum_{b\in
            B} r_t \cdot x_{tb}}$}
    && \hspace{3.9cm} (DM) \nonumber
    \\
    \makebox[13mm][l]{s.t.}
    \makebox[100mm][l]{$\displaystyle \sum_{b  \in B\setminus\{\beta\}}  a_{tb} \sum_{l\in L} P_l \cdot z_{b l}
-   a_{t \beta} \sum_{l\in L}  P_l \cdot z_{\beta l} +  M \cdot x_{t\beta} \leq
\delta + M
\hspace{1.8cm} t \in T, \hspace{0.1cm}
       \beta \in B
    $}
    \label{eq:DM-SIR}
    \\
    \makebox[100mm][l]{$\displaystyle{\sum_{b\in B} x_{tb} \leq 1}
    \hspace{9.1cm} t\in T
    $}
    \nonumber
    \\
    \makebox[100mm][l]{$\displaystyle \sum_{l \in L} z_{bl} = 1
    \hspace{9.15cm} b \in B
    $}
    \label{eq:DM-onepower}
    \\
    \makebox[100mm][l]{$x_{tb} \in \{0,1\}
    \hspace{9cm} t \in T, \hspace{0.1cm} b \in B
    $}
    \nonumber
    \\
    \makebox[100mm][l]{$z_{bl} \in \{0,1\}
    \hspace{9.05cm} b \in B, \hspace{0.1cm} l \in L.
    $}
    \label{eq:DM-z}
    \nonumber
\end{eqnarray}

\noindent Note that, due to (\ref{eq:powervalue}), every $p_b$
also satisfies $0\leq p_b \leq P_{max}$. As a consequence, the box
constraints on $p_b$ and thus variable $p_b$ are dropped from the
formulation.

The Power-Indexed formulation is obtained from (DM) by
substituting each \emph{knapsack SIR constraint} (\ref{eq:DM-SIR})
with a set of {\em GUB cover inequalities} \citep{Wo90}.

In the following, we denote a GUB cover inequality by the acronym \emph{GCI}. The GCIs constitute a stronger version of simple cover inequalities of a knapsack constraint, and are defined by exploiting the presence of the additional constraints (\ref{eq:DM-onepower}), that are called \emph{generalized upper bound} (GUB) constraints.

Before introducing the GCIs, we recall some related definitions and
concepts introduced in \citep{Wo90}. We consider the set of binary
points $Y = P \cap B^n$, where $P\subseteq R^n_+$ is the polytope
defined by:
\begin{eqnarray}
    \makebox[75mm][l]{($i$) $\displaystyle{\sum_{j\in N_1} a_j \cdot y_j - \sum_{j\in N_2} a_j \cdot  y_j \leq a_0}$}
     \nonumber
   \\[10pt]
    \makebox[25mm][l]{($ii$) $\displaystyle{\sum_{j\in S_i} y_j \leq 1 }$}
    \makebox[50mm][l]{ for  $i\in I_1 \cup I_2$}
    \label{eq:knapsack+gub}
    \\[10pt]
    \makebox[75mm][l]{$y\in R^n_+$,}
    \nonumber
 \end{eqnarray}

\noindent where $N = N_1 \cup N_2, \hspace{0.1cm} N_1\cap N_2 =\emptyset, \hspace{0.1cm} a_j
> 0$ for $j\in N$, $\bigcup_{i\in I_1} S_i = N_1$, $\bigcup_{i\in I_2}
S_i = N_2$ and, finally $S_i\cap S_l = \emptyset$ if $i,l \in I_k$
with $i\neq l$ for $k=1,2$.
In other words, the variables of the knapsack (\ref{eq:knapsack+gub}.i) are partitioned into
a number of subsets, and at most one variable can be set to 1 for each subset. Each of these subsets thus defines a GUB constraint (\ref{eq:knapsack+gub}.ii).  Furthermore, by definition each subset is entirely contained either in $N_1$ or $N_2$ and thus the coefficients of the corresponding 0-1 variables have the same sign in the knapsack constraint (\ref{eq:knapsack+gub}.i).


\noindent
A set $C = C_1 \cup C_2$ is a {\em GUB cover} for $Y$ if:
\begin{eqnarray}
    \makebox[100mm][l]{($i$) $\displaystyle{C_k \subseteq N_k}$
    \qquad for $k=1,2$}
    \nonumber
    \\[6pt]
    \makebox[100mm][l]{($ii$) $\displaystyle{|C_k\cap S_i| \leq 1}$ \hspace{0.2cm}
    for $i\in I_k$ and  $k=1,2$}
    \nonumber
    \\[6pt]
    \makebox[100mm][l]{($iii$) $\displaystyle{\sum_{j\in C_1} a_j - \sum_{j\in C_2} a_j > a_0.}$}
    \nonumber
 \end{eqnarray}

\noindent
On the basis of the GUB cover $C$, it is easy to build a standard cover inequality which is valid for the set $Y$.
Such constraint can be lifted by including new variables, by exploiting the GUB inequalities (\ref{eq:knapsack+gub}.ii). In particular,  with the GUB cover $C$ we associate the following sets:
\begin{eqnarray}
    \makebox[66mm][l]{$\displaystyle{I_k^+ = \{i\in I_k: C_k \cap S_i \neq \emptyset\}}$}
    \makebox[34mm][l]{for  $k=1,2$}
    \nonumber
    \\[6pt]
    \makebox[50mm][l]{$\displaystyle{S_i^+ = \{j\in S_i: a_j \geq a_l}$  for
    $l\in C_1 \cap S_i\}$}
    \makebox[50mm] {for $i\in I^+_1$}\nonumber \\[6pt]
    \makebox[50mm][l]{$\displaystyle{S_i^+ = \{j\in S_i: a_j \leq a_l}$  for
    $l\in C_2 \cap S_i\}$}
    \makebox[50mm] {for $i\in I^+_2.$}\nonumber
 \end{eqnarray}

\noindent
For each set $S_i$ with one element in the cover, $S_i^+$ represents the set of elements which may be added to the cover in order to lift the corresponding inequality. In particular, if the elements of $S_i$ correspond to non-negative coefficients $a_j$ of the knapsack, then we can add all the elements that correspond to coefficients that are larger than $a_l$ (i.e., the coefficient of the element of $S_i$ in the GUB cover). We instead include all the elements with smaller coefficient in the case of negative coefficients.

In \citep{Wo90}, Wolsey proves that if $C = C_1 \cup C_2$ is a {\em GUB cover}, the following GUB cover inequality (GCI) is valid for $Y$:
\begin{equation}\label{eq:Wolsey_gub-coverIneq}
    \sum_{i\in I_1^+} \sum_{j\in S_i^+} y_j \leq |C_1| - 1
    + \sum_{i\in I_2^+} \sum_{j \not\in S_i^+} y_j
    + \sum_{i\in I_2 \backslash I_2^+} \sum_{j \in S_i} y_j.
\end{equation}
\noindent
When $I_2^+ = I_2$ and $|I_2| = 1$, such valid inequality reduces to:
\begin{equation}\label{eq:gub-covers}
\sum_{i\in I_1^+} \sum_{j\in S_i^+}  y_j + \sum_{i\in I_2^+}
\sum_{j\in S_i^+} y_j \leq |C_1|.
\end{equation}

\medskip

\noindent
Now, let us focus on a single knapsack constraint
(\ref{eq:DM-SIR}) of (DM) associated with testpoint $t\in T$ and
server $\beta\in B$, along with  constraints
(\ref{eq:DM-onepower}) for $b \in B$ and the valid inequality
$x_{t\beta} \leq 1$. We can cast this into the GUB framework
introduced by Wolsey by making the following associations:
\begin{eqnarray}
    \makebox[100mm][l]{$\displaystyle{N_1 = \{(b,l): \hspace{0.1cm} b\in
    B \backslash \{\beta\}, \hspace{0.1cm} l\in L\}}
    \hspace{0.1cm} \cup \hspace{0.1cm} \{(t,\beta)\}$}
    \nonumber \\
    \makebox[100mm][l]{$\displaystyle{N_2 = \{(\beta,l): \hspace{0.1cm} l\in L.\}}$}
    \nonumber
 \end{eqnarray}

\noindent
Observe that, with a slight abuse of notation, in the
definition of $N_1$ we are also including index $(t,\beta)$
corresponding to variable $x_{t\beta}$. Similarly, we let:
\begin{eqnarray}
    \makebox[100mm][l]{$\displaystyle{I_1 = \{b: \hspace{0.1cm} b\in
    B \backslash \{\beta\}\}
    \hspace{0.1cm} \cup \hspace{0.1cm}
    \{(t,\beta)\}}$}
    \nonumber
    \\
    \makebox[100mm][l]{$\displaystyle{I_2 = \{\beta\}.}$} \nonumber
 \end{eqnarray}

\noindent
Indeed, for each $b\in B$ at most one variable
$z_{bl}$ can be equal 1, for $l\in L$, and we have $S_b = \{(b,l):
l\in L\}$ for all $b\in B$. Also, we let $S_{t,\beta} =
\{(t,\beta)\}$ be the singleton corresponding to variable
$x_{t\beta}$. Observe that we have $N_1 = S_{t,\beta} \cup (
\bigcup_{b \in B \backslash \{\beta\}} S_b )$ and $N_2 = S_{\beta}$.

Before translating conditions $(i)$, $(ii)$ and $(iii)$ into our
setting, we provide an intuitive explanation of how a CGI is build for formulation (DM). For a fixed couple of receiver and server and a fixed subset of interferers, a GUB cover corresponds to one serving power level and a combination of interfering power levels that jointly deny the coverage of the receiver by the server. Thereafter, the lifting is done by considering lower serving power levels and higher interfering power levels. We now proceed to define formally the GCI.
To this purpose, consider first the coverage condition
(\ref{eq:SIR}) corresponding to  receiver $t \in T$ with server
$\beta \in B$. Suppose that the server $\beta$ is emitting at
power value $p_\beta = P_{\lambda}$, for some $\lambda \in L$. Let
$\Gamma = \{b_1, \dots, b_{|\Gamma|}\} \subseteq
B\backslash\{\beta\}$ be a set of interferers (for $t$ when
$\beta$ is its server) and let $q_1, \dots, q_{|\Gamma|}$ be
power levels for each interferer in $\Gamma$ such that:
\begin{equation}\label{eq:not-served}
a_{t b_1} \cdot P_{q_1} + \dots  + a_{t b_{|\Gamma|}} \cdot P_{q_{|\Gamma|}} -
a_{t\beta} \cdot P_\lambda > \delta.
\end{equation}

\noindent
In other words, receiver $t$ is not served when $t$ is assigned to
server $\beta$ emitting at power value $P_{\lambda}$, and the
interferers $b_1, \dots, b_{|\Gamma|}$ are emitting at power
values $p_{b_1} = P_{q_1},$ $\dots,$ $p_{b_{|\Gamma|}} =
P_{q_{|\Gamma|}}$, respectively.

By letting $C_1 = \{(b_i,q_i): i = 1, \dots, |\Gamma|\}
\cup \{(t,\beta)\}$ and $C_2 = \{(\beta,\lambda)\}$, it follows that
$C = C_1 \cup C_2$ is a cover of (\ref{eq:DM-SIR}). Also, it is
not difficult to see that $C$ is a GUB cover, since $C_1 \subseteq
N_1$, $C_2 \subseteq N_2$, $|C_1 \cap S_b| \leq 1$, for all $b\in
I_1$ and $|C_2 \cap S_{\beta}| = 1$. We also have $I_1^+ = \Gamma
\cup \{(t,\beta)\}$ and $I_2^+ = \{\beta\}$.

\bigskip

\noindent
Since $a_{t b} \cdot P_{l} < a_{t b} \cdot P_{l+1}$ for all $b\in B$ and $l = 1, \dots, |L|-1$, we have that
 $S_{b_i}^+ = \{(b_i,q_i), (b_i, q_{i+1}), \dots, (b_i,q_{|L|})\}$ for $b_i \in
B \backslash \{\beta\}$, $S_{t,\beta}^+ = \{(t,\beta)\}$ and $S^+_{\beta}
= \{(\beta,1), \dots, (\beta,\lambda)\}$.

It follows from (\ref{eq:gub-covers}) that, for  $t\in T$, $\beta
\in B$, the inequality

\begin{equation}\label{eq:single-lifted-gub}
x_{t \beta} + \sum_{l = 1}^\lambda  z_{\beta l} +
\sum_{i=1}^{|\Gamma|} \sum_{j=q_i}^{|L|} z_{b_i j} \leq |\Gamma| +
1
\end{equation}

\medskip
\noindent
is valid for the set of binary vectors satisfying
(\ref{eq:DM-SIR}) and (\ref{eq:DM-onepower}).

\vspace{1cm}

\noindent
Now, for all the subsets of interferers $\Gamma \subseteq B \backslash
\{\beta\}$, denote by $L^I(t,\beta, \lambda, \Gamma)$ the set of
$|\Gamma|$-tuples $q \in L^{|\Gamma|}$ satisfying
(\ref{eq:not-served}). The following proposition follows
immediately by the validity of (\ref{eq:single-lifted-gub}):

\begin{proposition}\label{pro:DM-lifted-gub}
Given $t\in T$, $\beta \in B$, the family of inequalities:
\begin{equation} \label{eq:DM-lifted-gub}
x_{t \beta} + \sum_{l = 1}^\lambda  z_{\beta l} +
\sum_{i=1}^{|\Gamma|} \sum_{j=q_i}^{|L|} z_{b_i j} \leq |\Gamma| +
1
\end{equation}
defined for $\Gamma \subseteq
        B \backslash\{\beta\}$, $\lambda \in L$, $q \in L^I(t,\beta, \lambda,
        \Gamma)$,
is satisfied by all the binary solutions of (\ref{eq:DM-SIR}) and
(\ref{eq:DM-onepower}).

\end{proposition}

\noindent
It can be formally shown that the reverse is also true, namely all
binary solutions to (\ref{eq:DM-lifted-gub}) and
(\ref{eq:DM-onepower}) also satisfy (\ref{eq:DM-SIR}). It follows that the following formulation, that we
call Power-Indexed $(PI)$, is valid for the WND (with finite set
of power values):
\begin{eqnarray}
    \makebox[13mm][l]{$\max $}
    \makebox[120mm][l]{$\displaystyle{\sum_{t\in T} \sum_{b\in
            B} r_t \cdot x_{tb}}$}
    \nonumber
    && \hspace{2.1cm} (PI)
    \\
    \makebox[13mm][l]{s.t.}
    \makebox[120mm][l]{$\displaystyle{x_{t \beta} + \sum_{l = 1}^\lambda
    z_{\beta l} + \sum_{i=1}^{|\Gamma|} \sum_{j=q_i}^{|L|}
    z_{b_i j} \leq |\Gamma| + 1}
    \hspace{1.25cm}
    t\in T, \hspace{0.1cm}
    \beta \in B, \hspace{0.1cm}
    \Gamma \subseteq B \backslash \{\beta\},
    $}
    \nonumber
    \\
    \makebox[55mm][l]{$
    \hspace{0.50cm}
    \lambda \in L, \hspace{0.1cm} q \in L^I(t,\beta, \lambda, \Gamma)$}
    \label{eq:PI-multi-interferer}
    \\
    \makebox[65mm][l]{$\displaystyle{\sum_{b\in B} x_{tb} \leq 1}$}
    \makebox[55mm][l]{$\hspace{0.55cm} t \in T$}
    \label{eq:PI-oneserver}
    \\
    \makebox[65mm][l]{$\displaystyle \sum_{l \in L} z_{bl} = 1$}
    \makebox[55mm][l]{$\hspace{0.55cm} b \in B$}
    \label{eq:PI-onepower}
    \\
    \makebox[65mm][l]{$x_{tb} \in \{0,1\}$}
    \makebox[55mm][l]{$\hspace{0.55cm} t \in T, \hspace{0.1cm} b \in B$}
    \label{eq:PI-x}
    \\
    \makebox[65mm][l]{$z_{bl} \in \{0,1\}$}
    \makebox[55mm][l]{$\hspace{0.55cm} b \in B, \hspace{0.1cm} l \in L.$}
    \label{eq:PI-z}
\end{eqnarray}

\noindent
The above formulation contains a very large number of
GCIs (potentially exponential in $|B|$ for all $t\in T$). To
cope with this we proceed in a standard fashion by initially
considering a subset of all inequalities and subsequently
generating new inequalities when needed. In Section
\ref{sec:algorithm}, we give the details of our column and row
generation approach to solve the WND along with a heuristic
routine for separating violated GCIs
(\ref{eq:PI-multi-interferer}). The overall behaviour of the row
generation approach is strongly affected by the quality of the
initial relaxation. In the context of WND, a particularly
well-suited choice consists of including only the GCIs
(\ref{eq:PI-multi-interferer}) corresponding to interferer sets
$\Gamma$ with $|\Gamma| = 1$; we denote such initial relaxation by
$(PI^0)$. This choice has several major advantages.

First, the number of constraints in (PI$^0$) is small and can be
generated efficiently. In the next section, we actually show that,
for each $t\in T, \beta \in B$ and $b \in
B\backslash\{\beta\}$, the number of non-dominated GCIs (\ref{eq:PI-multi-interferer}) is at most $|L|$.

Second, as the Power-Indexed formulation (PI) is derived from
the discretized SIR formulation (DM), so (PI$^0$) can be thought
as derived from a relaxation (DM$^0$) of (DM). Namely, the
relaxation (DM$^0$) is obtained from (DM) by replacing, for each
$t\in T$ and each $\beta \in B$, the SIR inequality
(\ref{eq:DM-SIR}) with the family of inequalities (one for each
interferer):
\begin{equation}\label{eq:DM-binarySIR}
a_{tb} \sum_{l\in L} P_l \cdot z_{b l} -
a_{t \beta} \sum_{l\in L}  P_l \cdot z_{\beta l} + M
\cdot x_{t\beta} \leq \delta + M \qquad b\in B \backslash
\{\beta\}.
\end{equation}

\noindent Clearly, each inequality of type (\ref{eq:DM-binarySIR})
is dominated by the original inequality (\ref{eq:DM-SIR}) from
which it derives, and the 0-1 solutions to (DM$^0$) may not be
feasible for (DM). Nevertheless, in many applicative contexts
(DM$^0$) appears to be a very good approximation of (DM).
Indeed, this type of relaxation has been introduced in
\citep{MaMaSa09} to cope with DVB network design problems, and
successfully applied to the design of the Italian national
reference DVB network. Similarly, our experiments reported in Section
\ref{sec:results} show that (PI$^0$) is a good approximation of (PI).
Indeed, the number of inequalities not in (PI$^0$)
generated by our Branch-and-Cut is always very small.  This can be
well explained by the practical observation that,
for a given receiver, there exists most of the time one particular
interferer whose signal is much stronger than the others (see
 Section \ref{sec:results} for a more detailed discussion).

A third and most crucial feature of (PI$^0$)  relates to the strength
of its GCIs. In the next section we show that, for
each $t\in T$, $\beta \in B$ and $b \in
B \backslash\{\beta\}$, the family of GCIs associated with
(\ref{eq:DM-binarySIR}) along with the trivial facets define the
corresponding {\em GUB knapsack polytope}, i.e. the convex hull of
the 0-1 solutions to the knapsack SIR constraint
(\ref{eq:DM-binarySIR}) and its corresponding GUB constraints
(\ref{eq:DM-onepower}). This is a very desired property which
explains why the LP-relaxations of (PI$^0$) provide much tighter
bounds than those provided by (DM$^0$), thus in turn implying more
effective searches and the capability to solve larger instances.

Summarizing, (PI$^0$) can be easily generated, is a good
approximation of the original problem and provides strong
LP-relaxations.

\section{The GUB knapsack polytope for the single-interferer SIR inequality}  \label{sec:convex-hull}

For a receiver $t\in T$, server $\beta \in B$ and a single
interferer $b \in B\backslash\{\beta\}$, let us consider the
family of GCIs associated with the constraint
(\ref{eq:DM-binarySIR}):
\begin{equation} \label{eq:single-interf-gub-covers}
x_{t \beta} + \sum_{l = 1}^\lambda  z_{\beta l} +
 \sum_{j=q}^{|L|} z_{b j} \leq 2 \hspace{1.2cm}
 \lambda \in L, \hspace{0.1cm} q \in L^I(t,\beta, \lambda, \{b\}).
\end{equation}

\noindent Since $P_l > P_{l-1}$ for $q= 2,\dots, |L|$, the set
$L^I(t,\beta, \lambda, \{b\})$ of interfering levels of
$b$ for a server power level $\lambda$ can be written as
$\{q(\lambda), q(\lambda)+1, \dots, |L|\}$, where $q(\lambda) =
\min \{l \in L: a_{tb} \cdot  P_l - a_{t\beta}
\cdot P_\lambda
> \delta\} $.
It follows that the subfamily of inequalities
(\ref{eq:single-interf-gub-covers}) associated with $\lambda$ is
dominated by the single inequality corresponding to $q(\lambda)$.
Finally, observe that $q(\lambda') \geq q(\lambda)$ for $\lambda'
\geq \lambda$.

 In order to simplify the notation, we now let $u =
x_{t\beta}$, $v_l = z_{\beta l}$ for $l\in L$ and $w_l = z_{b
l}$ for $l \in L$. After removing the dominated GCIs, the
remaining family can be rewritten as:
\begin{equation} \label{eq:non-dominated-gub-covers}
u + \sum_{l = 1}^\lambda  v_l +
 \sum_{l=q(\lambda)}^{|L|} w_l \leq 2  \hspace{1.2cm}  \lambda = 1,\dots, |L|.
\end{equation}

\noindent
The following theorem extends a result presented in \citep{Wo90}
(Proposition 3.1), also providing an alternative and simpler proof
for it.

\begin{proposition}\label{th:convex-hull}
The polytope $P$ defined as the set of points $(u,{\bf v},{\bf w})
\in \mathbb{R}^{1+2|L|}$ satisfying (\ref{eq:non-dominated-gub-covers}) and
the  constraints $0\leq u\leq 1$, ${\bf 0}\leq {\bf v} \leq {\bf
1}$ and  ${\bf 0}\leq {\bf w} \leq {\bf 1}$ is the convex hull of
the 0-1 solutions to (\ref{eq:DM-binarySIR}).
\end{proposition}

\proof{Proof of Proposition \ref{th:convex-hull}.}
Let $A$ be the 0-1 coefficient matrix associated with the set of
constraints (\ref{eq:non-dominated-gub-covers}). We first show
that $A$ is an \emph{interval matrix}, i.e. in each column the 1's
appear consecutively \citep{NeWo88}.

We start by noticing that $A= (U|V|W)$ where $U$ is the column
associated with the  variable $u$; $V \in \{0,1\}^{|L|\times|L|}$
is the square matrix associated with the variables $v_1, \dots,
v_{|L|}$; and  $W \in \{0,1\}^{|L|\times|L|}$  is the square
matrix associated with the  variables $w_1, \dots, w_{|L|}$.

The vector $U$ has all the elements equal to 1 as $u$ is included
in every constraint (\ref{eq:non-dominated-gub-covers}). The
matrix $V = [n_{ij}]$ with $i,j = 1,\ldots,|L|$ is lower
triangular and such that $n_{ij} = 1$ for $i \geq j$.
Indeed, the constraint (\ref{eq:non-dominated-gub-covers})
corresponding with $\lambda \in L$ includes exactly the $v$
variables $v_1, \dots, v_{\lambda}$.

Finally, consider the matrix $W = [m_{ij}]$ with $i,j =
1,\ldots,|L|$. First, observe that for all $\lambda, j \in L$, we
have:
$$
m_{\lambda j} = 1 \hspace{0.1cm} \Longleftrightarrow
\hspace{0.1cm} j \geq q(\lambda).
$$

Recalling that for every $\lambda', \lambda \in L$ with $\lambda'
\geq \lambda$, we have $ q(\lambda') \geq q(\lambda)$, it follows
that,  for all $\lambda \leq \lambda'$,  $m_{\lambda' j} = 1
\Longrightarrow j \geq q(\lambda') \Longrightarrow j \geq
q(\lambda) \Longrightarrow m_{\lambda j} = 1$. $W$ is thus an
interval matrix and as $U$ and $V$ are interval matrices as well,
it follows that $A$ is an interval matrix and thus totally
unimodular.

Finally, if we denote by $\bar{A}$ the matrix associated with the
constraints (\ref{eq:non-dominated-gub-covers}) and the box
constraints on variables $u, {\bf v}, {\bf w}$, then $\bar{A}$ is
obtained by extending $A$ with $I$ and $-I$, where $I$ is the
identity matrix of size $1+2|L|$. Thus $\bar{A}$ is a totally unimodular
matrix \citep{NeWo88} and, since the right hand sides of the
constraints are integral, the vertices of $P$ are also integral,
completing the proof.\Halmos
\endproof

\section{Solution Algorithm} \label{sec:algorithm}

The solution algorithm is based on the (PI) formulation for the WND
and consists of two basic steps: $(i)$ a set $\cal P$ of feasible
power values is established; $(ii)$ the associated formulation is
solved by row generation and Branch-and-Cut. We start by
describing step $(ii)$ and we come back to step (i) later in this
section.

In the following, for a fixed power set ${\cal P}$, we denote the
solution algorithm for the associated (PI) formulation as
SOLVE-PI(${\cal P}$).  Since the (PI) formulation has in general
an exponential number of constraints of type
(\ref{eq:PI-multi-interferer}), we apply row generation. Namely,
we start by considering only a suitable subset of constraints and
we solve  the associated relaxation. We then check if any of the
neglected rows is violated by the current fractional solution. If
so, we add the violated row to the formulation and solve again,
otherwise we proceed with  standard Branch-and-Cut (as implemented
by the commercial solver Cplex). The separation
of violated constraints is repeated in each branching node.

At node 0, the initial formulation (PI$^0$) includes only a subset
of constraints (\ref{eq:PI-multi-interferer}), namely those including one
interferer (i.e., $|\Gamma| = 1$). In Section
\ref{sec:power-indexed} and Section \ref{sec:convex-hull} we
discussed why this is a good choice for (PI$^0$). Indeed, in our case
studies, only a low number of additional constraints is added by separation during the iterations of the algorithm.

\subsection{Separation.}

We now proceed to show how violated constraints are separated. Let
$(x^*,z^*)$ be the current fractional solution. In Section
\ref{sec:power-indexed} we have showed that constraints
(\ref{eq:PI-multi-interferer}) are GUB cover inequalities
of (\ref{eq:DM-SIR}). In order to separate a violated GCI of type (\ref{eq:PI-multi-interferer}), we make use of the exact oracle introduced by  \cite{Wo90} and heuristically solve it by extending the
standard (heuristic) approach to the separation of cover
inequalities described in \citep{NeWo88}.

To this end, let us first select a receiver $t \in T$ and one of
its servers, say $\beta \in B$. We want to find a GCI of type (\ref{eq:PI-multi-interferer}) that is
associated with $t$ and $\beta$,  and is violated by the current
solution $(x^*,z^*)$. In other words, we want to identify a power
level $\lambda \in L$ for $\beta$, a set of interferers $\Gamma =
\{b_1, \dots, b_{|\Gamma|}\}\subseteq B \backslash \{\beta\}$
and an interfering $|\Gamma|$-tuple of power levels $q = (q_1,
\dots, q_{|\Gamma|}) \in L^I(t,\beta, \lambda, \Gamma)$, such
that:
\begin{equation}\label{eq:violated}
 {x^*_{t\beta} + \sum_{l = 1}^\lambda z^*_{\beta l} +
\sum_{i=1}^{|\Gamma|} \sum_{j=q_i}^{|L|} z^*_{b_i j} > |\Gamma| +
1}.
\end{equation}

\noindent
Recall that $q \in L^I(t,\beta,
\lambda, \Gamma)$ if
\begin{equation}\label{eq:q-interferer}
\sum_{i = 1}^{|\Gamma|} a_{t b_i} \cdot P_{q_i} - a_{t \beta} \cdot P_\lambda > \delta.
\end{equation}

\noindent
We solve the separation problem by defining a suitable
0-1 Linear Program. In particular, in order to identify a suitable
pair $(\beta,\lambda)$ we introduce, for every $l\in L$, a
binary variable $u_{\beta l}$, which is $1$ iff $l = \lambda$.
Similarly, we introduce binary variables $u_{bl}$ for all $b\in
B \backslash \{\beta\}$ and $l\in L$, with $u_{bl} = 1$ iff
$(b,l)= (b_i,q_i)$, where $b_i\in \Gamma$ and $q_i$ is the
corresponding interfering power level. Then $u\in
\{0,1\}^{|B(t)|\times |L|}$ satisfies the following system of
linear inequalities:
\begin{equation}\label{eq:cover-oracle}
\sum_{b\in B \backslash \{\beta\}} a_{tb} \hspace{0.1cm} \sum_{l\in L}  P_{l} \cdot u_{bl} -
a_{t \beta} \hspace{0.1cm} \sum_{l \in L} P_l \cdot u_{\beta l}  > \delta
\end{equation}
\begin{equation}\label{eq:onepower-oracle}
\sum_{l\in L}  u_{bl} = 1 \qquad b\in B.
\end{equation}

\noindent
Constraint (\ref{eq:cover-oracle}) ensures that $u$ is
the incidence vector of a cover of (\ref{eq:DM-SIR}), whereas
constraint (\ref{eq:onepower-oracle}) states that $u$ satisfies
the GUB constraints.

Observe now that $|\Gamma| = \sum_{b \in B \setminus \{\beta\}}
\sum_{l\in L} u_{bl}$. So,  if $u$ identifies a violated GCI (\ref{eq:violated}), we must have:
\begin{equation}\label{eq:fo-oracle}
\sum_{l \in L} u_{\beta l} \sum_{k = 1}^l z^*_{\beta k} +
\sum_{b\in B \backslash \{\beta\}} \sum_{l \in L} u_{b l}
\sum_{k=l}^{|L|} z^*_{b k}
\hspace{0.2cm} > \hspace{0.2cm}
\sum_{b\in B \backslash \{\beta\}}
\sum_{l\in L}  u_{b l} + 1 - x^*_{t \beta}.
\end{equation}

\medskip

\noindent
In order to (heuristically) search for a violated
inequality, we proceed in a way which resembles the classical
approach for standard cover inequalities \citep{NeWo88}, by
considering the following linear program (SEP), introduced by \cite{Wo90}:
\begin{eqnarray}
 Z &=&
 \max
 \hspace{0.5cm}
 \sum_{l\in L} u_{\beta l} \sum_{k =
    1}^l z^*_{\beta k} + \sum_{b\in B \setminus \{\beta\}}
    \sum_{l\in L} u_{b l} \cdot
    \left(\sum_{k=l}^{|{L}|} z^*_{b k}-1\right)
\nonumber
 \hspace{2.7cm} (SEP)
\\
& &
\nonumber
\\
& & \mbox{s.t.}
\hspace{0.2cm} \sum_{b\in B \setminus \{\beta\}} a_{tb} \sum_{l\in L}
P_{l} \cdot u_{bl} - a_{t \beta} \sum_{l \in L} P_{l} \cdot u_{\beta l} \geq \delta
\label{eq:knapsack-SEP}
\\
& &
\hspace{0.8cm} \sum_{l\in L} u_{bl} = 1
\hspace{6.4cm} b\in B
\nonumber
\\
& &
\hspace{0.8cm}
u_{bl} \geq 0
\hspace{7.0cm}
b\in B, \hspace{0.1cm} l\in L.
\nonumber
\end{eqnarray}

\noindent
It is easy to notice that the feasible region of (SEP)
contains all binary vectors satisfying (\ref{eq:cover-oracle}) and
(\ref{eq:onepower-oracle}). Let $Z$ be the optimum value to (SEP).
If $Z \leq 1 - x^*_{t \beta}$ then no binary vector $u$ satisfies
(\ref{eq:fo-oracle}) and consequently no violated constraint
exists. If $Z > 1 - x^*_{t \beta}$ then a violated constraint may
exist, and we resort to a heuristic approach to find it. In
particular, observe first that $Z$ can be computed by relaxing the
knapsack constraint (\ref{eq:knapsack-SEP}) in a Lagrangian
fashion and then by solving the resulting Lagrangian dual, namely:
\begin{eqnarray*}
Z & = & \min_{\eta \geq 0} \hspace{0.3cm} Z(\eta)
\end{eqnarray*}

\noindent
where $\eta\in \mathbb{R}^+$ is the Lagrangian multiplier
and:
\begin{eqnarray*}
Z(\eta) & = & \max_{u \geq 0} \hspace{0.3cm} \sum_{l\in L}
u_{\beta l} \sum_{k = 1}^l z^*_{\beta k} + \sum_{b\in B \setminus
\{\beta\}} \sum_{l\in L} u_{b l} \cdot \left(\sum_{k=l}^{|{L}|} z^*_{b
k}-1\right)
\\
& & + \hspace{0.3cm} \eta \cdot \left( \sum_{b\in B \setminus \{\beta\}} a_{tb} \sum_{l\in L}
P_{l} \cdot u_{bl} - a_{t\beta} \sum_{l \in L} P_{l} \cdot u_{\beta l} - \delta \right)
\\
& &
\\
& & \mbox{s.t.} \hspace{0.4cm} \sum_{l\in L} u_{bl} = 1 \qquad
b\in B.
\end{eqnarray*}

\medskip

\noindent For fixed $\eta \geq 0$, the objective $Z(\eta)$ can be
easily computed by inspection. To simplify the notation we rewrite
the objective function of the above linear program as:
\begin{equation}\label{eq:objsep}
-\delta \cdot \eta \hspace{0.2cm} + \hspace{0.2cm} \max_{u\geq 0}
\hspace{0.1cm} \sum_{b\in B} \sum_{l\in L} c_{bl}(\eta) \cdot
u_{bl}
\end{equation}

\noindent
where, for every $b \in B, l \in L$, we let:
\begin{displaymath}
c_{bl}(\eta)= \left\{
\begin{array}{ll}
\sum_{k = 1}^l z^*_{\beta k} -
\eta \cdot a_{t \beta} \cdot P_{l}  & \hspace{0.7cm} \mbox{if  } b = \beta
\\
\\
\sum_{k=l}^{|L|} z^*_{b k} - 1  + \eta \cdot a_{tb} \cdot P_{l} &
\hspace{0.7cm} \mbox {if  } b\in B \backslash \{\beta\}.
\\
\end{array}
\right.
\end{displaymath}

\medskip

\noindent
For fixed $\eta \geq 0$, an optimal solution $u(\eta)$
to the inner maximization problem can be found by inspection as
follows. For each $b\in B$, identify a power level $l_b \in L$
which maximizes the coefficient in (\ref{eq:objsep}), namely
$c_{bl_b}(\eta) = \max_{l \in L} c_{bl}(\eta)$; then,
 for each $b\in B$ and each $l\in L$, let:
\begin{equation}\nonumber
 u_{bl}(\eta) = \left\{
\begin{array} {cl}
1 & \mbox{\hspace{0.2cm} if $l = l_b$}\\
0 & \mbox{\hspace{0.2cm} otherwise}.
\end{array}
\right.
\end{equation}

\smallskip

\noindent It is straightforward to see that, for all $\eta \geq
0$, $u(\eta) \geq 0$ satisfies all constraints
(\ref{eq:onepower-oracle}) and maximizes  (\ref{eq:objsep}). For
$\eta \geq 0$, the function $Z(\eta)$ is convex and unimodal and the optimum solution $\eta^*$ can be found
efficiently by applying the {\em Golden Section Search Method}
\citep{GeWh04}. Suppose now that $Z(\eta^*) > 1 - x^*_{t
\beta}$ (otherwise no violated constraints exist). If, in
addition,  $u(\eta^*)$ also satisfies (\ref{eq:cover-oracle}),
 then the positive components of the
binary solution $u(\eta^*)$ are in one-to-one correspondence to
the variables of a violated constraint. Otherwise the algorithm
returns no violated cover.

Finally, when the current solution $(x^*,z^*)$ is purely 0-1, we perform an exact separation by directly checking the satisfaction of each of the constraints (\ref{eq:PI-multi-interferer}).

\bigskip

\subsection{The Algorithm}
We come back now to the first step in our algorithm, namely the
choice of the set of admissible power values $\cal P$. Large sets
are in principle more likely to produce better quality solutions.
However, the ability of the solution algorithm to find optimal or
simply good-quality solutions is strongly affected by $|{\cal
P}|$, as we will show in more details in the computational results
section. Thus, the size and the elements of $\cal P$ should
represent a suitable compromise between these two opposite
behaviors. Moreover, the effectiveness of the Branch-and-Cut is
typically affected by the availability of a good initial feasible
solution. Thus, we decided to iteratively apply SOLVE-PI(${\cal
P}$) to a sequence of  power sets ${\cal P}_0 \subset {\cal P}_1
\subset \dots \subset {\cal P}_r$. Each invocation inherits all
the generated cuts, the best solution found so far and the
corresponding lower bound from the previous invocation. More
precisely, if we denote by -99 the switched-off state  (in dBm),
and $P_{\min}^{dBm}$, $P_{\max}^{dBm}$ are the (integer) minimum
and maximum power values (in dBm), then we have ${\cal P}_0 =
\{-99, P_{\max}^{dBm}\}$, ${\cal P}_1 = \{-99,
P_{\min}^{dBm},\left\lfloor\frac
{P_{\max}^{dBm}-P_{\min}^{dBm}}{2}\right\rfloor, P_{\max}^{dBm}\}$
and ${\cal P}_r = \{-99, P_{\min}^{dBm}, P_{\min}^{dBm} + 1,
\dots, P_{\max}^{dBm}\}$. The structure of the intermediate power
sets will be described in Section \ref{sec:results}. Observe that
the actual power values are only used in the separation oracle
where the dB values are converted into the original non-dB values.

The overall approach, denominated \emph{WPLAN}, is summarized in
Algorithm \ref{wplan}, where $i$ denotes the current iteration,
along with the associated best solution found $x_i$, the
corresponding value $LB_i$, and the set of feasible powers  ${\cal
P}_i$.  If SOLVE-PI(${\cal P}_i$) is executed in less than the
iteration time limit TL$_{i}$ then the residual time $\tau_i$ is
used to increase the time limit of the following iteration (i.e.,
TL$_{i+1}$:= TL$_{i+1}$ + $\tau_i$). The initial incumbent
solution $x_{-1}$ corresponds to all transmitters switched off and
no receiver served ($LB_{-1} = 0$).
\renewcommand{\algorithmicrequire}{\textbf{Input:}}
\renewcommand{\algorithmicensure}{\textbf{Output:}}
\begin{algorithm}
\caption{\mbox{WPLAN}} \label{wplan}
\begin{algorithmic}
\vspace{0.2cm} \REQUIRE the power sets  ${\cal P}_0, {\cal P}_1,
\dots, {\cal P}_r$
, the iteration time limit TL$_{i} \mbox{ for } i=0,\ldots, r$
\\
\ENSURE the best solution $x_r$


\STATE $LB_{-1} := 0$
 \FOR{$i=0$ to $r$}
        \STATE \textbf{1. } Invoke SOLVE-PI(${\cal P}_i$) with lower bound $LB_{i-1}$, incumbent $x_{i-1}$ and TL$_{i}$
        \STATE \textbf{2. } Get $x_i$, $LB_{i}$ and $\tau_{i}$
            \STATE \textbf{3. } $TL_{i+1} := TL_{i+1} + \tau_i$
 \ENDFOR
 \STATE Return $x_r$
\end{algorithmic}
\end{algorithm}

\section{Computational Results} \label{sec:results}

The model that we have considered so far has a very simple and basic structure and applies to the main wireless technologies. More precisely, it can be effectively used if the service coverage condition of a receiver is expressed by means of a SIR constraint (\ref{eq:firstSIRinequality}). As pointed out in Section \ref{sec:power-indexed}, each technology generally requires its own peculiar parameter values and additional constraints and/or variables to model its own specific features.

In this section, we present computational results concerning realistic instances of two important wireless technologies: the \emph{IEEE Standard 802.16} \citep{WiMAXstd04} and the Terrestrial Digital Video Broadcasting (DVB-T, \citealt{ETSI06}).

The target of these tests is manyfold. First, we compare the new
(PI) formulation to the two big-\emph{M} formulations (BM) and
(DM) and show that (PI) outperforms (BM) and (DM) both in terms
of quality of bounds and quality of solutions. Then, we
illustrate specific features of the solution algorithm WPLAN and
we motivate the iterative approach with increasing power sets.
Finally,  we assess the ability of WPLAN to tackle realistic
network design instances. The tests were performed under Windows
XP 5.1 operating system, with 1.80 GHz Intel Core 2 Duo processor
and 2$\times$1024 MB DDR2-SD RAM. The algorithm is implemented in
C++ (under Microsoft Visual Studio 2005 8.0), whereas the
commercial MILP solver ILOG Cplex 10.1 is invoked by ILOG Concert
Technology 2.3.

In the following two subsections, we provide a concise description of the main specific features of the two technologies and we highlight their impact on the basic model that we presented in Section \ref{sec:power-indexed}. Furthermore, we describe the characteristics of the realistic instances that we consider for each technology.

\subsubsection{WiMAX Network Design}

The first set of instances refers to a \emph{WiMAX network} and were developed with the \emph{Technical Strategy and Innovations Unit} of \emph{British Telecom Italia} (BT). WiMAX is the common name used to indicate the \emph{IEEE Standard 802.16} \citep{WiMAXstd04}. Specifically, we consider the design of a \emph{Fixed WiMAX Network} that provides broadband internet access.}

The major amendments concern the introduction of different
frequency channels, channel capacity and traffic demand. To model
the additional features of a WiMAX network, the formulations (BM)
and (PI) must include additional variables to take into account
multiple frequencies (denoted by set $F$) and multiple
transmission schemes (denoted by set $H$). Furthermore, we need to
introduce  additional constraints to model the capacity of each
frequency to accommodate traffic generated by users.  For a
detailed description of these additional features, both from
technological and modeling perspective, we refer the reader to
\citep{Zh09}.
%
%
%

All the instances correspond to an urban area of the city of Rome (Italy), selected in agreement with the engineers at BT,
who considered it as a representative residential traffic
scenario. Each activated transmitter can emit by using integer power levels in the range [20,40] dBm. We define three types of instances, denoted by SX, where X is the instance identifier ranging in $\{1,\dots,7\}$, RX with X
$= \{1,\dots,4\}$ and QX with X $= \{1,\dots,4\}$. For the SX
instances, the traffic is uniformly distributed among the TPs and
we assign unitary revenue to each TP (i.e. $r_t =1$). Finding an
optimal coverage plan thus corresponds to define the plan with the
maximum number of covered TPs. Only one frequency and one burst
profile are allowed. For the RX instances,  we consider a traffic
distribution based on the actual distribution of the buildings. We
also introduce multiple frequencies and burst profiles. In this
case, the revenue of each testpoint is proportional to the traffic
generated. Finally, the QX instances include an increasing number
of candidate sites and focus on a single frequency network with
multiple burst profiles. The dimension of each instance is resumed
in Table \ref{tab:WiMAXinstances}.

\begin{table}[htbp]
\caption{Description of the WiMAX test-bed instances}
\label{tab:WiMAXinstances} \small
\begin{center}
\begin{tabular}{|c||ccccccc|cccc|cccc|}
  \hline
  ID & S1 & S2 & S3 & S4 & S5 & S6 & S7 & R1 & R2 & R3 & R4 & Q1 & Q2 & Q3 & Q4
  \\
  \hline
  \hline
  $|\mbox{T}|$ & 100 & 169 & 196 & 225 & 289 & 361 & 400 & 400 & 441 & 484 & 529 & 400 & 441 & 484 & 529
  \\
  $|\mbox{B}|$ & 12 & 12 & 12 & 12 & 12 & 12 & 18 & 18 & 18 & 27 & 27 & 36 & 36 & 36 & 36
  \\
  $|\mbox{F}|$ & 1 & 1 & 1 & 1 & 1 & 1 & 1 & 3 & 3 & 3 & 3 & 1 & 1 & 1 & 1
  \\
  $|\mbox{H}|$ & 1 & 1 & 1 & 1 & 1 & 1 & 1 & 4 & 4 & 4 & 4 & 4 & 4 & 4 & 4
  \\
  \hline
\end{tabular}
\end{center}
\end{table}

\subsubsection{DVB-T Network Design}

The second set of instances refers to networks based on the \emph{Terrestrial Digital Video Broadcasting} technology (DVB-T, \citealt{ETSI06}). Indeed, our algorithm has been used to design the \emph{reference networks of the Italian DVB-T plan}, comprising 25 national and hundreds of regional single-frequency networks. Unfortunately, due to non-disclosure agreements, we cannot reproduce and distribute the details of the real life instances. Nevertheless, we have synthesized 9 instances using the same digital terrain and propagation model, the same population database and, finally, the same technical assumptions defined by the \emph{Italian Authority for Telecommunications} (Agcom). As a consequence, our instances and solutions constitute a valid proxy of the real networks planned by the Authority and currently under deployment by the italian broadcasters.

Each instance corresponds to a regional area of Italy, with an extent ranging from about 3.500 to about 30.000 km$^2$.  The network represented in an instance is constituted by a set of transmitters $B$ that synchronously broadcast the same telecommunication service on the same frequency over a target area. Each transmitter can emit by using a subset of power levels in the range [-40,26] dBkW. Service coverage is evaluated in a set of testpoints $T$ and the revenue obtained by covering a testpoint is equal to the population living in the corresponding elementary portion of territory.
The coverage is assessed through an adapted version of the SIR inequality (\ref{eq:SIR}): the rules of distinction between serving and interfering signals and summation of signals comes from the adoption of \emph{Orthogonal Frequency Division Multiplexing} (OFDM) in the DVB-T technology. For a detailed description of how the SIR inequality is built, we refer the reader to \cite{MaRoSm06}. The dimension of each instance is shown in Table \ref{tab:DVBinstances}.

We stress that the coefficient matrices associated to the DVB-T instances are in general more ill-conditioned than those associated to the WiMAX instances. This can be intuitively explained by considering that DVB-T networks involve transmitters that are much more powerful than those used by a WiMAX network. Such transmitters are able to broadcast signals at very long distance. As a consequence, weak signals can be picked up also far away from the target area, creating interference that may be very small when compared to the powerful signal of closer serving transmitters (for example, Italian transmitters in Sardegna may interfere transmissions in Tunisia and Southern France). The ratio between the largest and the smallest fading coefficient of a DVB-T SIR inequality is thus in general much larger than that of a WiMAX SIR inequality. Numerical instability phenomena become therefore more marked.

\begin{table}[htbp]
\caption{Description of the DVB-T test-bed instances}
\label{tab:DVBinstances} \small
\begin{center}
\begin{tabular}{|c||ccccccccc|}
  \hline
  ID & DVB1 & DVB2 & DVB3 & DVB4 & DVB5 & DVB6 & DVB7 & DVB8 & DVB9
  \\
  \hline
  \hline
  $|\mbox{T}|$   & 2003 & 1741 & 5618 & 4466 & 2704 & 4421 & 197 & 3400 & 2003
  \\
  $|\mbox{B}|$   &  127 &  188 &  411 &  202 &  113 &  215 & 109 &  183 & 127
  \\
  \hline
\end{tabular}
\end{center}
\end{table}

\subsection{Numerical Results and Comparisons}
\label{subsec:results}

We have pointed out in Section \ref{sec:introduction} that the
solutions to (BM) and (DM) returned by state-of-the-art MILP
solvers such as Cplex can be affected by numerical inaccuracy,
i.e. the SIR inequalities of testpoints recognized as covered are
actually unsatisfied (similar problems were also reported in
\cite{KaKeOl06}, \cite{KeOlRa10} and \cite{MaMaSa09}). We detect such coverage
errors by evaluating the solutions \emph{off-line}: after the
optimization process, we verify that the SIR inequality
corresponding to each nominally covered testpoint is really
satisfied by the power vector of the returned solution. This is
not the only issue, as, in the case of some instances, (BM) and (DM) can be even wrongly evaluated as infeasible.

In our experience, tuning the parameters of Cplex is crucial to
reduce coverage errors and to contain the effects of numerical
instability. Furthermore, in the case of (DM), tuning is
essential to ensure that the problem is correctly recognized as
feasible. After a series of tests, we established that, in the case of (BM) and (DM), an
effective setting consists of turning off the \emph{presolve} and
on the \emph{numerical emphasis}. Moreover, we turn off the
generation of the \emph{mixed-integer rounding cuts} and of the
\emph{Gomory fractional cuts} as we observed no advantages in the
quality of the bounds and a sensible increase in running times.

\subsubsection*{Assessing the strength of the Power-Indexed formulation.}
\label{subsection:strengthPI}

The first group of experiments is designed
to assess the strength of (PI)  comparing it with (BM) and
(DM). To this end, we focus on a single instance of our test-bed
(instance S4 presented in Table \ref{tab:WiMAXinstances}) and
detail the behaviour of WPLAN for each invocation of
 SOLVE-PI(${\cal P}$). The sets of power values in the first three invocations of
SOLVE-PI(${\cal P}$) are (in dBm) ${\cal P}_1 = \{-99,40\}, {\cal
P}_2 = \{-99,20, 30, 40\}$ and ${\cal P}_3 = \{-99, 20, 25, 30, 35,
40\}$, respectively. Then, in each of the following invocations,
${\cal P}$ is expanded  by including two more values (suitably
spaced). To analyse the behaviour of the single iterations and
establish an effective sequence of power sets, we set a time limit
of 1 hour for each invocation of the solution algorithm for (PI)
and (DM).

In order to evaluate the quality of (PI) w.r.t. (DM), we apply
WPLAN to (DM) (note that in this case the solution procedure
SOLVE-PI is replaced by the simple solution of (DM) by
Cplex). In Table \ref{tab:behaviour}, for each iteration of WPLAN,
we report the number $|\mbox{L}|$ of considered power levels, the
number of GCIs included in the initial formulation (PI$^0$) and
the number of GCIs separated during the current iteration.
Additionally, for both (PI) and (DM), we report the upper
bound at node 0 $(UB)$, the value $|\mbox{T*}|$ of the final
solution (number of covered testpoints) and the final gap. When
the solution contains coverage errors, two values are presented in
the $|\mbox{T*}|$ column, namely the nominal value of the best
solution returned by Cplex (in brackets) and its actual value
computed by re-evaluating the solution \emph{off-line}.

The last line of the table shows the results obtained for (BM)
by setting a time limit of 3 hours. Note that in this case, the
second column reports the number of SIR (big-\emph{M}) constraints
(\ref{eq:BM-SIR}) included in (BM). This number is by definition also the number of SIR (big-\emph{M}) constraints (\ref{eq:DM-SIR}) included in (DM).

\begin{table}[htbp]
\caption{Behaviour of WPLAN for instance S4} \label{tab:behaviour}
\begin{center}
\begin{tabular}{|c||cc|ccc||ccc|}
  \hline
        &   \multicolumn{2}{c|}{GCIs}    &
\multicolumn{3}{c||}{(PI)} & \multicolumn{3}{c|}{(DM)}
    \\
      $|\mbox{L}|$ & init &  added     &  UB  &  $|\mbox{T*}|$  & gap\%
&  UB  &  $|\mbox{T*}|$  & gap\%
\\
   \hline
   \hline
   2    &   5743   &  17 &   199.2193    &   106  &   0.00
   &   218.3465    &   91  &  125.65
\\
        4    &  9035   &  7 &   204.2500    &     111  &   0.00
&   219.0015    &   97 (98)  &   102.68
   \\
        6    &  14312   &  13 &   206.6261    &     111  &   59.03
&   219.3488    &   100 (101)  &   115.70
   \\
        8    &  17142   &  45   &   209.4200    &   111  &   67.51
&   219.7349    &   100 (101)  &   122.98
   \\
       10    &  24638   &  6  &   210.0000    &    111   &   79.99
&   220.2788    &   100 (101)  &   123.14
   \\
        12   &  27799   &  1  &   211.7000    &    111  & 82.05
&   219.9144    &   100 (101)  &   124.01
   \\
       14   &  35944   &  0  &   212.0000    &     111    &   83.46
&   220.1307    &   100 (101)  &   123.58
   \\
       16   &  38496   &  10  &   214.5930    &     111    &   85.48
&   220.3000    &  100 (101)  &   125.00
   \\
        18   &  45425   &  2  &   215.8000    &    111    &   86.44
&   220.1091    &   100 (101)  &   124.83
   \\
        20   &  48918   &  2  &   218.0000    &    111    &   89.99
&   220.0560    &  100 (101)  &   125.00
   \\
        22   &  57753   &  3  &   218.0000    &    111    &   90.83
&   220.3720    &  100 (101)  &   125.00
   \\
  \hline
  \hline
        (BM)        &   1170   &  -  &   221.3925    &    93   &
 97.18 &   -    &   -  &   -
 \\
  \hline
\end{tabular}
\end{center}
\end{table}

The figures in Table \ref{tab:behaviour} are representative of the
typical behaviour of WPLAN on all instances of our test-bed. They
allow us to make some relevant observations. First, the size of
(PI) grows quickly with the number of power levels, and is
typically much larger than that of (BM) and (DM). This is
counterbalanced by the quality of the upper bounds, which are
consistently better for (PI) and, most important, the quality of
the solutions found. Interestingly, the best solution is found
quite early in the iterative procedure, namely for $|{\cal P}|
\leq 6$. A similar behaviour is observed for the other WiMAX instances
reported in Table  \ref{tab:comparisonsWiMAX} and the DVB-T instances in Table \ref{tab:comparisonsDVB} as well. This motivated
our choice of the sequence of feasible power values in the final
version of WPLAN for WiMAX: most of the computational effort is concentrated on small cardinality power sets, and only one large set. More precisely, there will be only 4 iterations, corresponding to 2, 4, 6 and 22 power levels, respectively.

Finally, we note that the number of generated GCIs is small.
Also, in most cases the GCIs include only two interferers, and in
any case never more than three. In other words, even though many
interferers can reach a given testpoint, only very few of them (in
most cases only one) give a significant contribution to the
overall interference.

\subsubsection*{The performance of the Power-Indexed approach over the test-bed.}

In this subsection, we comment the results over our WiMAX and DVB-T benchmark instances. For an exhaustive report of the results through tables, we refer the reader to the Appendix of this paper. In all experiments, we set a time limit of 3 hours for the solution of (BM) and (DM) and for WPLAN applied to (PI).
As in subsection \ref{subsection:strengthPI},
we solve (DM) by an adapted version of WPLAN (we recall that in this case the solution procedure SOLVE-PI is replaced by the simple solution of (DM) by Cplex).

Besides the results obtained by solving the ``pure'' models (BM) and (DM), we report also the results obtained by trying to stabilize (BM) and (DM) through Cplex \emph{indicator constraints} and by strengthening (DM) through a suitable subset of our GCIs.
The indicator constraints constitute a way to express relationships between variables and may reduce the flaws of big-{\em M} formulations. In our work, we check if
declaring the big-{\em M} coverage constraints of (BM) and (DM) by means of Cplex indicator constraints \citep{CPLEX} can improve the quality of solutions. We denote the resulting formulations by adding the symbol ``+'' to the acronym (e.g., BM+).
Furthermore, we investigate if it is convenient to strengthen (DM) by simply including the GCIs corresponding with the single-interferer condition (i.e., $|\Gamma| = 1$, see Section \ref{sec:algorithm}). We denote the resulting formulation by \emph{(DM \& GCI1)}. This investigation is motivated by the fact that such subset of GCIs seems to be very effective to discover high quality solutions fast.

The results show that WPLAN applied to (PI) outperforms (BM)
and (DM) in terms of quality of the solutions found and,  in most cases, running times to obtain them (the running times obviously include also the
time spent by the separation oracle). Coverage errors, in particular, are completely eliminated. Even if in principle the
reduced and quite small number of power values considered by WPLAN
could result in poorer coverage w.r.t. (BM), the results clearly
show that this is not the case. On one hand, this happens as a
small number of well-spaced power values suffices in practice to
obtain good coverage; indeed, it is common practice in WiMAX
network planning to neglect intermediate values, i.e. a device is
either switched-off or activated at its maximum power
\citep{Ri10}. On the other hand, the size of the (BM)
formulation and the ill-conditioned constraint matrix, along with
the presence of the big-\emph{M} coefficients, makes the solution
process unstable, the solutions found unreliable and the
branching tree extremely large. Indeed, due to rounding errors and
numerical instability, several solutions to (BM) turn out to be
infeasible when verified off-line. The effects of numerical instability become more marked in the case of the DVB-T test bed: the feasible solutions to (BM) of all but one of the instances contain coverage errors that entail the loss of up to 20.000 users.
Furthermore, in contrast to the good performance of WPLAN, several instances seem to be very difficult for (BM) and no feasible solution is retrieved by the time limit.

WPLAN applied to (PI) also outperforms (DM) solved by the adapted WPLAN algorithm. The results show that in general the simple
discretization of the power range does not suffice to get better
solutions than those obtained by (BM). Indeed, in many cases, the performance of (DM) is worse than that of (BM) and coverage errors are still strongly present.


Finally, after having pointed out the advantages of a pure GCI-based approach, we assess if stabilizing (DM) by indicator constraints or strengthening (DM) by GCIs can lead to remarkable advantages.
The results reported show that stabilizing by indicator constraints does not allow to reach the quality of the solutions obtainable by the pure GCI formulation (PI).
%
This  behavior can be explained by the simple observations that (PI) allows us to get rid of the major sources of instability in (BM) and (DM), namely  the bad-conditioned coefficients of the constraint matrix (not affected by the use of Cplex indicator constraints) and the big-\emph{M} coefficients.
%
Though in a significative number of cases the value of the best solution is higher than that of solutions obtained by pure (BM) and (DM), coverage errors are still (heavily) present. Moreover, the value of the solutions is anyway lower than those obtained by (PI).
Strengthening (DM) by GCIs seems to be more effective than stabilizing by indicator constraints:
in many cases (DM \& GCI1) reaches better solutions w.r.t. (DM) and (DM+) However, also in this case, there are still a few solutions that contain errors and the final value is anyway lower than that obtained through (PI). Finally, also in the case of (DM+ \& GCI1), the stabilization by indicator constraints seems to decrease the value of solutions obtained within the time limit, without being able to completely avoid coverage errors.

\subsubsection*{Comparisons between warm and cold start for (PI).}

Finally, in Table \ref{tab:warmcold} we show the impact of the
iterative approach WPLAN on the quality of the solutions found for
$(PI)$ in the case of the WiMAX instances. A similar behaviour is observed also in the case of the DVB-T instances.
In particular we compare {\em cold starts}, which
correspond to invoking SOLVE-PI(${\cal P}$) without benefiting
from cuts and lower bounds obtained at former invocations, with
{\em warm starts} which, in contrast, make use of such
information. The value of the best solutions found during
successive invocations of SOLVE-PI both under warm and cold starts
are shown in the columns identified by $|L| = n$, where $n$
denotes the number of corresponding power levels. The value of the
best solution found at the first invocation is in column $|L| =
2$, while the value of the best solution and the number of levels
used to find it are shown in column $|T^*|$ and $|L^*|$,
respectively.

\begin{table}[htbp]
\caption{Comparisons between warm and cold starts}
\label{tab:warmcold}
\begin{center}
\begin{tabular}{|c|ccc|cc|cc|}
  \hline
  & & & &\multicolumn{2}{c|}{WARM START} & \multicolumn{2}{c|}{COLD START}
  \\
  \raisebox{1.5ex}{ID} & \raisebox{1.5ex}{$|\mbox{T*}|$} & \raisebox{1.5ex}{$|\mbox{L*}|$} & \raisebox{1.5ex}{$|\mbox{L}|$=2} & $|\mbox{L}|$=4 & $|\mbox{L}|$=6 & $|\mbox{L}|$=4 & $|\mbox{L}|$=6
  \\
  \hline
  \hline
  S1  &  74 & 6 &  69 &  72 &  74 &  71 &  58 \\
  S2  & 107 & 4 &  72 & 107 & 107 &  80 &  63 \\
  S3  & 113 & 4 &  83 & 113 & 113 & 108 &  101 \\
  S4  & 111 & 4 &  75 & 111 & 111 & 100 &  97 \\
  S5  &  86 & 6 &  76 &  84 &  86 &  83 &  81 \\
  S6  & 170 & 4 & 127 & 170 & 170 & 110 & 127 \\
  S7  & 341 & 4 & 296 & 341 & 341 & 314 & 196  \\
  \hline
  R1  & 400 & 2 & 400 & - & - & 399 & 304\\
  R2  & 441 & 4 & 416 & 441 & - & 394 & 355 \\
  R3  & 427 & 2 & 427 & 427 & 427 & 414 & Out\\
  R4  & 529 & 2 & 529 & - & - & 512 & Out\\
  \hline
  Q1  & 67 & 2 & 67 & 67 & 67 & $\ast$ & $\ast$ \\
  Q2  & 211 & 4 & 196 & 211 & 211 & 156 & Out  \\
  Q3  & 463 & 2 & 463 & 463 & 463 & Out & Out  \\
  Q4  & 491 & 2 & 491 & 491 & 491 & Out & Out \\
  \hline
\end{tabular}
\end{center}
\end{table}

For all S-instances the best solution can be found only due to
warm start. Note that SOLVE-PI encounters increasing difficulties
in finding good solutions as the number of power levels increases
(in the case of the apparently hard instance Q1, for 3 and 5 power
levels, no feasible solution is found within the time limit when
cold start is adopted). This is mainly due to the large size of
the corresponding instances, that, in some cases denoted by
\emph{Out}, makes Cplex run out of memory while building the
model.  However, a good initial solution provided to SOLVE-PI can
be improved in most cases. We have already observed that for a
larger number of levels  (i.e. $ > 6$), no improved solutions can
be found for all the WiMAX instances in our test-bed. Finally, for R1 and R4 a solution covering the entire target area is found already with
$|L| = 2$, while for R2 such a solution is found with $|L|=4$ (and
warm-start).

\section{Conclusions.}
The coverage condition in wireless network design problems is typically modeled by linearizing the signal-to-interference ratio and by including the notorious big-M coefficients. The resulting Mixed-Integer Programs are very weak and ill-conditioned, hence unable to solve large instances of real networks.
In this paper, we show how power discretization, a common modeling approach among professionals, can constitute the first step to define formulations that are noticeably stronger than the classical ones. These pure 0-1 formulations are based on GUB cover inequalities, that completely eliminates the source of numerical instability. This new Power-Indexed approach outperforms the
classical big-{\em M}  models, both in terms of quality of solutions found and of strength of the bound, as showed by an extensive computational study on real WiMAX and DVB-T instances.


%
%
%

\begin{APPENDIX}{Tables of comparison.}
\label{appendix}

In this appendix, we present tables that exhaustively report the computational results about comparisons between (BM), (DM) and (PI),  that we have commented in Section \ref{subsec:results}.
The results are reported in Tables \ref{tab:comparisonsWiMAX} and \ref{tab:comparisonsWiMAX_indicCnstr} for the WiMAX instances and in Tables \ref{tab:comparisonsDVB} and \ref{tab:comparisonsDVB_indicCnstr} for the DVB-T instances.
All results are obtained by setting a time limit of 3 hours.
We recall that (BM+), (DM+) respectively denote
the versions of (BM),(DM) stabilized through Cplex \emph{indicator constraints}, whereas (DM \& GCI1) denotes the version of (DM) strengthened through the GCIs corresponding with the single-interferer condition (i.e., $|\Gamma| = 1$).

The value of the best solutions found within the time limit is shown in column $|\mbox{T*}|$ for WiMAX and column \emph{COV\%} for DVB-T (\emph{COV\%} is the percentage of population covered with service). The {\em gap\%} columns report the \emph{nominal} (i.e., before checking solution correctness) percentage gap between the upper and lower bound at termination, the {\em time} column specifies when the best solution is found (in \emph{seconds}), whereas the last column $|\mbox{L*}|$ is the number of power levels used in the iteration in which WPLAN obtains the best solution. In the columns $|\mbox{T*}|$ and \emph{COV\%}, the expression ``Out'' indicates that Cplex runs out of memory while building the model. Finally, to denote some specials situations, we adopt the following conventions: i) the expression \emph{``None by TL''} indicates that the solver is not able to find a feasible solution within the time limit; ii) the expression \emph{``Infeasible$^\ast$''} indicates that the solver wrongly considers the problem as being infeasible.

We briefly resume the three main observations that can be made on the basis of the results (discussed in more detail in Section \ref{subsec:results}): 1) WPLAN applied to (PI) outperforms (BM) and (DM) for all the WiMAX and DVB-T instances; 2) in most cases, the use of indicator constraint leads to finding solutions of lower value than those provided by pure (DM) and this reduction in value is not compensated by a complete elimination of coverage errors; 3) strengthening (DM) by GCIs in general enhances the solving performance, but solutions containing errors are still generated.

The higher performance of our approach based on the Power-Indexed formulation (PI) is especially apparent for all of the DVB instances and the WiMAX
R-instances. In particular, several instances seem to be quite easy for WPLAN but very difficult for (BM) and (DM). Indeed, when no time limit is imposed to the solution of (BM), Cplex runs out of memory after
about ten hours of computation without getting sensible
improvements in the bounds. On the contrary, in the case of WiMAX instances like R1, R2
and R4, SOLVE-PI(${\cal P}$) finds the optimum solution (when
$|{\cal P}| = 2$) in less than 1 hour. The higher performance is
also highlighted in the case of instance Q1 that turns  out to be
hard: both (BM) and (DM) with 2 power levels cannot find any
feasible solution with non-zero value within the time limit,
while, in contrast, (PI) finds a solution with value 67.

\begin{table}[htbp]
\caption{Comparisons between (BM), (DM), (DM \& GCI1) and  WPLAN (WiMAX instances)}
\label{tab:comparisonsWiMAX} \small
\begin{center}
\begin{tabular}{|cc||ccc|ccc|ccc|ccc|}
  \hline
  & &
  \multicolumn{3}{c|}{(BM)}
  & \multicolumn{3}{c|}{(DM)}
  & \multicolumn{3}{c|}{(DM \& GCI1)}
  & \multicolumn{3}{c|}{WPLAN}
  \\
  \raisebox{1.5ex}{ID} & \raisebox{1.5ex}{$|\mbox{T}|$} &
  $|\mbox{T*}|$ & gap\% (nom) & time &
  $|\mbox{T*}|$ &  time & $|\mbox{L*}|$ &
  $|\mbox{T*}|$ &  time & $|\mbox{L*}|$ &
  $|\mbox{T*}|$ &  time & $|\mbox{L*}|$
  \\
  \hline
  \hline
  S1 & 100 & 63 (78)  &  13.72 & 10698 &
  65 (70) & 8705  & 6 &
  67 & 8275 & 6 &
  74  & 10565 & 6
  \\
  S2 & 169 & 99 (100) &  56.18 & 10705  &
  86 & 6830  & 4 &
  101 & 6405 & 4 &
  107  & 5591 & 4
  \\
  S3 & 196 & 108 & 79.54 & 4010  &
  61 (100) & 4811 & 4 &
  95 & 6908 & 4 &
  113 & 5732 & 4
  \\
  S4 & 225 &  93  & 103.43 & 10761  &
  100 (101) & 6507 & 6 &
  100 & 6891 & 4 &
  111  &  7935  & 4
  \\
  S5 & 289 &  77  & 202.24 & 10002  &
  73 (76) & 4602 & 4 &
  72 (82) & 7088 & 4 &
  86  & 10329  & 6
  \\
  S6 & 361 & 154  & 130.76 &  8110  &
  121 (138) & 5310 & 4 &
  149 & 5724 & 4 &
  170  &  8723  & 4
  \\
  S7 & 400 & 259 (266) & 49.67 & 8860 &
  120 (121) & 4100 & 4 &
  239 & 6003 & 4 &
  341 & 7154 & 4
  \\
  \hline
  R1 & 400 & 370 &   7.57 & 10626  &
  284 (328) & 1066 & 2 &
  304 & 3424 & 2 &
  400  & 1579  & 2
  \\
  R2 & 441 & 302 (303)  &  45.03 &  3595  &
  393 (394) & 4713 & 4 &
  375 (384) & 4371 & 4 &
  441  & 1244  & 4
  \\
  R3 & 484 &  99  & 385.86 & 10757  &
  188 & 2891 & 2 &
  306 & 3440 & 2 &
  427  & 3472  & 2
  \\
  R4 & 529 & 283 (286)  &  84.96 & 10765 &
  307 & 3026 & 2 &
  399 & 3152 & 2 &
  529  & 2984  & 2
  \\
  \hline
  Q1 & 400 &  0  & - & - &
  0 & - & - &
  37 & 3108 & 2 &
    67 & 2756 & 2
  \\
  Q2 & 441 &  191  & 130.89 & 9124 &
  158 (179) & 6282 & 4 &
  156 & 6932 & 4 &
  211  & 7132 & 4
  \\
  Q3 & 484 &  226  & 112.83 & 3392 &
  290 (292) & 2307 & 2 &
  316 & 3091 & 2 &
    463 &  3323  & 2
  \\
  Q4 & 529 & 145 (147) & 264.83 & 6623 &
  273 (280) & 1409 & 2 &
  343 & 2248 & 2 &
   491 & 3053 & 2
  \\
  \hline
\end{tabular}
\end{center}
\end{table}

\begin{table}[htbp]
\caption{Comparisons between (BM+), (DM+), (DM+ \& GCI1) and WPLAN (WiMAX instances)}
\label{tab:comparisonsWiMAX_indicCnstr} \small
\begin{center}
\begin{tabular}{|cc||ccc|ccc|ccc|ccc|}
  \hline
  & &
  \multicolumn{3}{c|}{(BM+)}
  & \multicolumn{3}{c|}{(DM+)}
  & \multicolumn{3}{c|}{(DM+ \& GCI1)}
  & \multicolumn{3}{c|}{WPLAN}
  \\
  \raisebox{1.5ex}{ID} & \raisebox{1.5ex}{$|\mbox{T}|$} &
  $|\mbox{T*}|$ & gap\% (nom) & time &
  $|\mbox{T*}|$ &  time & $|\mbox{L*}|$ &
  $|\mbox{T*}|$ &  time & $|\mbox{L*}|$ &
  $|\mbox{T*}|$ &  time & $|\mbox{L*}|$
  \\
  \hline
  \hline
  S1 & 100 & 63 (71) & 24.48 & 7556 &
  53 (56) & 7551 & 6 &
  68 (70) & 8650 & 6 &
  74  & 10565  & 6
  \\
  S2 & 169 & 99 (100) & 65.92 & 10322 &
  73 & 7036 & 4 &
  67 & 8816 & 4 &
  107  &  5591  & 4
  \\
  S3 & 196 & 101 (103) & 85.43 & 10241 &
  Infeasible$^*$ & - & - &
  72 (75) & 7101 & 4 &
  113 & 5732 & 4
  \\
  S4 & 225 & 71 & 179.2 & 8213 &
  97 (102) & 6820 & 4 &
  93 & 6566 & 4 &
  111  &  7935  & 4
  \\
  S5 & 289 & 69 & 262.55 & 6561 &
  75 & 6695 & 6 &
  75 & 6710 & 4 &
  86  & 10329  & 6
  \\
  S6 & 361 & 81 (107) & 235.51 & 5630 &
  98 (116) & 5672 & 4 &
  108 & 6102 & 4 &
  170  &  8723  & 4
  \\
  S7 & 400 & 238 & 67.23 & 7141 &
  157 (241) & 6008 & 4 &
  158 & 7059 & 4 &
  341 & 7154 & 4
  \\
  \hline
  R1 & 400 & 309 (340) & 17.05 & 4320 &
  292 & 2355 & 2 &
  188 & 3004 & 2 &
  400  & 1579  & 2
  \\
  R2 & 441 & 329 & 30.06 & 7809 &
  296 & 4682 & 4 &
  371 & 4798 & 4 &
  441  & 1244  & 4
  \\
  R3 & 484 & 0 & - & - &
  185 (201) & 3150 & 2 &
  374 & 3256 & 2 &
  427  & 3472  & 2
  \\
  R4 & 529 & 249 (253) & 112.44 & 9203 &
  278 & 2933 & 2 &
  329 & 2871 & 2 &
  529  & 2984  & 2
  \\
  \hline
  Q1 & 400 & 0 & - & - &
  0 & - & - &
  53 & 3400 & 2 &
    67 & 2756 & 2
  \\
  Q2 & 441 & 115 & 283.47 & 6135 &
  137 & 5024 & 4 &
  128 & 7082 & 4 &
  211  & 7132 & 4
  \\
  Q3 & 484 & 238 (263) & 82.88 & 5162 &
  258 & 2644 & 2 &
  291 & 2808 & 2 &
    463 &  3323  & 2
  \\
  Q4 & 529 & 252 & 109.92 & 7212 &
  378 (416) & 2118 & 2 &
  341 & 2391 & 2 &
   491 & 3053 & 2
  \\
  \hline
\end{tabular}
\end{center}
\end{table}

\begin{table}[htbp]
\caption{Comparisons between (BM), (DM), (DM \& GCI1) and WPLAN (DVB-T instances)}
\label{tab:comparisonsDVB} \small
\begin{center}
\begin{tabular}{|c||ccc|ccc|ccc|ccc|}
  \hline
  & \multicolumn{3}{c|}{(BM)}
  & \multicolumn{3}{c|}{(DM)}
  & \multicolumn{3}{c|}{(DM \& GCI1)}
  & \multicolumn{3}{c|}{WPLAN}
  \\
  \raisebox{1.5ex}{ID} &
  COV\% & gap\% & time  &
  COV\% &  time  & $|\mbox{L*}|$ &
  COV\% &  time  & $|\mbox{L*}|$ &
  COV\% &  time  & $|\mbox{L*}|$
  \\
  \hline
  \hline
  DVB1 &
  95.10 (95.49) & 2.67 & 8707 &
  96.76 & 10020 & 6 &
  94.80 & 7817 & 6 &
  97.26 & 7193 & 6
  \\
  DVB2 &
  96.15 (96.51) & 1.35 & 9510 &
  89.75 (96.98) & 8639 & 6 &
  96.03 & 9194 & 6 &
  97.14 & 9305 & 6
  \\
  DVB3 & None by TL & - & - &
  70.43 & 6155 & 4 &
  70.80 & 6891 & 4 &
  71.08  & 6544 & 4
  \\
  DVB4 & None by TL & - & - &
  77.55 (83.50) & 5650 & 4 &
  80.19 (84.85) & 7123 & 4 &
  88.98 & 725 & 2
  \\
  DVB5 & 94.94 (96.18) & 0.68 & 9804 &
  92.25 (94.30) & 7701 & 6 &
  95.93 (96.11) & 8826 & 8 &
  96.25  & 9677 & 8
  \\
  DVB6 & None by TL & - & - &
  71.73 & 5922 & 4 &
  74.09 & 6421 & 4 &
  74.55  & 5840 & 4
  \\
  DVB7 & 94.82 (100.00) & 0.00 & 65 &
  80.49 (100.00) & 293 & 8 &
  96.84 (99.83) & 311 & 10 &
  99.47 & 182 & 9
  \\
  DVB8 & 78.63 & 22.24 & 6621 &
  84.35 & 6212 & 4 &
  83.54 & 6049 & 4 &
  84.35 & 6293 & 4
  \\
  DVB9 & None by TL & - & - &
  95.26 & 5003 & 4 &
  95.57 & 5447 & 4 &
  96.60 & 1538 & 2
  \\
  \hline
\end{tabular}
\end{center}
\end{table}

\begin{table}[htbp]
\caption{Comparisons between (BM+), (DM+), (DM+ \& GCI1) and WPLAN (DVB-T instances)}
\label{tab:comparisonsDVB_indicCnstr}
\small
\begin{center}
\begin{tabular}{|c||ccc|ccc|ccc|ccc|}
  \hline
  & \multicolumn{3}{c|}{(BM+)}
  & \multicolumn{3}{c|}{(DM+)}
  & \multicolumn{3}{c|}{(DM+ \& GCI1)}
  & \multicolumn{3}{c|}{WPLAN}
  \\
  \raisebox{1.5ex}{ID} &
  COV\% & gap\% & time  &
  COV\% &  time  & $|\mbox{L*}|$ &
  COV\% &  time  & $|\mbox{L*}|$ &
  COV\% &  time  & $|\mbox{L*}|$
  \\
  \hline
  \hline
  DVB1 &
  95.38 (95.50) & 2.85 & 9203 &
  93.45 (96.62) & 9642 & 6 &
  94.27 & 8559 & 6 &
  97.26 & 7193 & 6
  \\
  DVB2 &
  96.95 & 0.90 & 9940 &
  96.22 (96.77) & 10077 & 6 &
  94.59 & 9803 & 6 &
  97.14 & 9305 & 6
  \\
  DVB3 & None by TL & - & - &
  69.94 & 6774 & 4 &
  70.21 & 7105 & 4 &
  71.08  & 6544 & 4
  \\
  DVB4 & 65.53 (65.65)  & 49.18 & 9120 &
  83.07 (86.67) & 6910 & 4 &
  77.57 (86.02) & 7008 & 4 &
  88.98 & 725 & 2
  \\
  DVB5 & 95.14 (95.41) & 1.61 & 6194 &
  94.03 & 9180 & 6 &
  95.02 & 10092 & 8 &
  96.25  & 9677 & 8
  \\
  DVB6 & None by TL & - & - &
  70.56 (70.96) & 6701 & 4 &
  70.30 & 6338 & 4 &
  74.55  & 5840 & 4
  \\
  DVB7 & 96.91 (100.00) & 0.00 & 244 &
  92.87 (99.81) & 573 & 10 &
  98.10 & 414 & 10 &
  99.47 & 182 & 9
  \\
  DVB8 & 58.51 & 64.24 & 10086 &
  80.09 (80.12) & 7050 & 4 &
  83.80 & 7122 & 4 &
  84.35 & 6293 & 4
  \\
  DVB9 & None by TL & - & - &
  95.85 & 6108 & 4 &
  94.44 & 6966 & 4 &
  96.60 & 1538 & 2
  \\
  \hline
\end{tabular}
\end{center}
\end{table}

\end{APPENDIX}

\ACKNOWLEDGMENT{The authors gratefully acknowledge the existence of
the Journal of Irreproducible Results and the support of the Society
for the Preservation of Inane Research.}





\bibliographystyle{nonumber}

\end{document}